\documentclass{amsproc}
\usepackage{amssymb}
\usepackage{amsmath}
\usepackage{eurosym}
\usepackage{amsfonts}

            \advance \topmargin by -\headheight \advance
            \topmargin by -\headsep
            \evensidemargin
            \oddsidemargin
\newtheorem{theorem}{Theorem}[section]
\newtheorem{lemma}[theorem]{Lemma}
\newtheorem{corollary}[theorem]{Corollary}

\newtheorem{remark}[theorem]{Remark}
\numberwithin{equation}{section}
\input{tcilatex}

\begin{document}
\title[On Generalized Tribonacci Sedenions]{On Generalized Tribonacci
Sedenions}
\thanks{}
\author[Y\"{u}ksel~Soykan, \.{I}nci Okumu\c{s} and Erkan Ta\c{s}demir]{Y\"{u}%
ksel~Soykan, \.{I}nci Okumu\c{s} and Erkan Ta\c{s}demir}
\maketitle

\begin{center}
\textsl{Department of Mathematics, }

\textsl{Art and Science Faculty, }

\textsl{Zonguldak B\"{u}lent Ecevit University,}

\textsl{67100, Zonguldak, Turkey }

\textsl{e-mail: \ yuksel\_soykan@hotmail.com}

\textsl{inci\_okumus\_90@hotmail.com}

\textsl{erkantasdemir@hotmail.com}
\end{center}

\textbf{Abstract.} The sedenions form a 16-dimensional Cayley-Dickson
algebra. In this work, we introduce the generalized Tribonacci sedenion and
present some properties of this sedenion.



\section{Introduction}

The generalized Tribonacci sequence $\{V_{n}(V_{0},V_{1},V_{2};r,s,t)\}_{n%
\geq 0}$ (or shortly $\{V_{n}\}_{n\geq 0}$) is defined as follows:%
\begin{equation}
V_{n}=rV_{n-1}+sV_{n-2}+tV_{n-3},\text{ \ \ \ \ }V_{0}=a,V_{1}=b,V_{2}=c,%
\text{ \ }n\geq 3  \label{equation:nbvcftyopuybnaeuo}
\end{equation}%
where $V_{0},V_{1},V_{2}\ $are arbitrary integers and $r,s,t$ are real
numbers.

This sequence has been studied by many authors, see for example [\ref%
{bib:bruce1984}], [\ref{catalani2002}], [\ref{bib:choi2013}],[\ref%
{bib:elia2001}],[\ref{bib:lin1988}], [ \ref{bib:pethe1988}], [\ref%
{bib:scott1977}], [\ref{bib:shannon1977}], [\ref{bib:spickerman1981}], [\ref%
{bib:yalavigi1972}], [\ref{bib:yilmaz 2014}].

The sequence $\{V_{n}\}_{n\geq 0}$ can be extended to negative subscripts by
defining%
\begin{equation*}
V_{-n}=-\frac{s}{t}V_{-(n-1)}-\frac{r}{t}V_{-(n-2)}+\frac{1}{t}V_{-(n-3)}
\end{equation*}%
for $n=1,2,3,...$ when $t\neq 0.$ Therefore, recurrence (\ref%
{equation:nbvcftyopuybnaeuo}) holds for all integer $n.$

If we set $r=s=t=1$ and $V_{0}=0,V_{1}=1,V_{2}=1$ then $\{V_{n}\}$ is the
well-known Tribonacci sequence and if we set $r=s=t=1$ and $%
V_{0}=3,V_{1}=1,V_{2}=3$ then $\{V_{n}\}$ is the well-known Tribonacci-Lucas
sequence.

In fact, the generalized Tribonacci sequence is the generalization of the
well-known sequences like Tribonacci, Tribonacci-Lucas, Padovan
(Cordonnier), Perrin, Padovan-Perrin, Narayana, third order Jacobsthal and
third order Jacobsthal-Lucas.\ In literature, for example, the following
names and notations are used for the special case of $r,s$ and $t.$

\begin{equation*}
\begin{array}{cccc}
\text{Sequences (Numbers)} &  & \text{Notation} &  \\ 
\text{Tribonacci} &  & \{T_{n}\}=\{V_{n}(0,1,1;1,1,1)\} &  \\ 
\text{Tribonacci-Lucas} &  & \{K_{n}\}=\{V_{n}(3,1,3;1,1,1)\} &  \\ 
\text{Padovan (Cordonnier)} &  & \{P_{n}\}=\{V_{n}(1,1,1;0,1,1)\} &  \\ 
\text{Pell-Padovan} &  & \{V_{n}(1,1,1;0,2,1)\} &  \\ 
\text{Jacobsthal-Padovan} &  & \{V_{n}(1,1,1;0,1,2)\} &  \\ 
\text{Perrin} &  & \{R_{n}\}=\{V_{n}(3,0,2;0,1,1)\} &  \\ 
\text{Pell-Perrin} &  & \{V_{n}(3,0,2;0,2,1)\} &  \\ 
\text{Jacobsthal-Perrin} &  & \{V_{n}(3,0,2;0,1,2)\} &  \\ 
\text{Padovan-Perrin} &  & \{S_{n}\}=\{V_{n}(0,0,1;0,1,1)\} &  \\ 
\text{Narayana} &  & \{N_{n}\}=\{V_{n}(0,1,1;1,0,1)\} &  \\ 
\text{third order Jacobsthal} &  & \{J_{n}\}=\{V_{n}(0,1,1;1,1,2)\} &  \\ 
\text{third order Jacobsthal-Lucas} &  & \{j_{n}\}=\{V_{n}(2,1,5;1,1,2)\} & 
\end{array}%
\end{equation*}%
The first few values of the sequences with non-negative indices are shown
below.%
\begin{equation*}
\begin{array}{ccccccccccccc}
n & 0 & 1 & 2 & 3 & 4 & 5 & 6 & 7 & 8 & 9 & 10 & ... \\ 
T_{n} & 0 & 1 & 1 & 2 & 4 & 7 & 13 & 24 & 44 & 81 & 149 & ... \\ 
K_{n} & 3 & 1 & 3 & 7 & 11 & 21 & 39 & 71 & 131 & 241 & 443 & ... \\ 
P_{n} & 1 & 1 & 1 & 2 & 2 & 3 & 4 & 5 & 7 & 9 & 12 & ... \\ 
R_{n} & 3 & 0 & 2 & 3 & 2 & 5 & 5 & 7 & 10 & 12 & 17 & ... \\ 
S_{n} & 0 & 0 & 1 & 0 & 1 & 1 & 1 & 2 & 2 & 3 & 4 & ... \\ 
N_{n} & 0 & 1 & 1 & 1 & 2 & 3 & 4 & 6 & 9 & 13 & 19 & ... \\ 
J_{n} & 0 & 1 & 1 & 2 & 5 & 9 & 18 & 37 & 73 & 146 & 293 & ... \\ 
j_{n} & 2 & 1 & 5 & 10 & 17 & 37 & 74 & 145 & 293 & 586 & 1169 & ...%
\end{array}%
\end{equation*}%
The first few values of the sequences with negative indices are shown below.%
\begin{equation*}
\begin{array}{ccccccccccccc}
n & 0 & 1 & 2 & 3 & 4 & 5 & 6 & 7 & 8 & 9 & 10 & ... \\ 
T_{-n} & 0 & 0 & 1 & -1 & 0 & 2 & -3 & 1 & 4 & -8 & 5 & ... \\ 
K_{-n} & 3 & -1 & -1 & 5 & -5 & -1 & 11 & -15 & 3 & 23 & -41 & ... \\ 
P_{-n} & 1 & 0 & 1 & 0 & 0 & 1 & -1 & 1 & 0 & -1 & 2 & ... \\ 
R_{-n} & 3 & -1 & 1 & 2 & -3 & 4 & -2 & -1 & 5 & -7 & 6 & ...%
\end{array}%
\end{equation*}%
As $\{V_{n}\}$ is a third order recurrance sequence (difference equation),
it's characteristic equation is $x^{3}-rx^{2}-sx-t=0,$ whose roots are%
\begin{eqnarray*}
\alpha &=&\alpha (r,s,t)=\frac{r}{3}+A+B \\
\beta &=&\beta (r,s,t)=\frac{r}{3}+\omega A+\omega ^{2}B \\
\gamma &=&\gamma (r,s,t)=\frac{r}{3}+\omega ^{2}A+\omega B
\end{eqnarray*}%
where%
\begin{eqnarray*}
A &=&\left( \frac{r^{3}}{27}+\frac{rs}{6}+\frac{t}{2}+\sqrt{\Delta }\right)
^{1/3},\text{ }B=\left( \frac{r^{3}}{27}+\frac{rs}{6}+\frac{t}{2}-\sqrt{%
\Delta }\right) ^{1/3} \\
\Delta &=&\Delta (r,s,t)=\frac{r^{3}t}{27}-\frac{r^{2}s^{2}}{108}+\frac{rst}{%
6}-\frac{s^{3}}{27}+\frac{t^{2}}{4},\text{ \ }\omega =\frac{-1+i\sqrt{3}}{2}%
=\exp (2\pi i/3)
\end{eqnarray*}%
Note that we have the following identities%
\begin{eqnarray*}
\alpha +\beta +\gamma &=&r, \\
\alpha \beta +\alpha \gamma +\beta \gamma &=&-s, \\
\alpha \beta \gamma &=&t.
\end{eqnarray*}%
From now on, we assume that $\Delta (r,s,t)>0,$ so that the Equ. (\ref%
{equation:nbvcftyopuybnaeuo}) has one real ($\alpha $) and two non-real
solutions with the latter being conjugate complex. So, in this case, it is
well known that generalized Tribonacci numbers can be expressed, for all
integers $n,$ using Binet's formula%
\begin{equation}
V_{n}=\frac{P\alpha ^{n}}{(\alpha -\beta )(\alpha -\gamma )}+\frac{Q\beta
^{n}}{(\beta -\alpha )(\beta -\gamma )}+\frac{R\gamma ^{n}}{(\gamma -\alpha
)(\gamma -\beta )}  \label{equat:mnopcvbedcxzsa}
\end{equation}%
where%
\begin{equation*}
P=V_{2}-(\beta +\gamma )V_{1}+\beta \gamma V_{0},\text{ }Q=V_{2}-(\alpha
+\gamma )V_{1}+\alpha \gamma V_{0},\text{ }R=V_{2}-(\alpha +\beta
)V_{1}+\alpha \beta V_{0}.
\end{equation*}%
Note that the Binet form of a sequence satisfying (\ref%
{equation:nbvcftyopuybnaeuo}) for non-negative integers is valid for all
integers $n,$ for a proof of this result see [\ref{bib:howard2010}]$.$This
result of Howard and Saidak [\ref{bib:howard2010}] is even true in the case
of higher-order recurrence relations.

We can give Binet's formula of the generalized Tribonacci numbers for the
negative subscripts:\ for $n=1,2,3,...$ we have 
\begin{equation*}
V_{-n}=\frac{\allowbreak \alpha ^{2}-r\alpha -s}{t}\frac{P\alpha ^{1-n}}{%
(\alpha -\beta )(\alpha -\gamma )}+\frac{\allowbreak \beta ^{2}-r\beta -s}{t}%
\frac{Q\beta ^{1-n}}{(\beta -\alpha )(\beta -\gamma )}+\frac{\allowbreak
\gamma ^{2}-r\gamma -s}{t}\frac{R\gamma ^{1-n}}{(\gamma -\alpha )(\gamma
-\beta )}.
\end{equation*}

In this paper, we intruduce and study the generalized Tribonacci sedenions
in the next section and give some properties of them. Before giving it's
definition, we present some information on Cayley-Dickson algebras.

The process of extending the field $\mathbb{R}$ of reals to $\mathbb{C}$ is
known as Cayley--Dickson process (construction). It can be carried further
to higher dimensions using a doubling procedure, yielding the quaternions $%
\mathbb{H}$, octonions $\mathbb{O}$ and sedenions $\mathbb{S}$ which (as a
real vector space) are of dimension $4$, $8$, and $16$, respectively. This
doubling process can be extended beyond the sedenions to form what are known
as the $2^{n}$-ions (see\ for example [\ref{bib:imaeda2000}], [\ref%
{bib:moreno1998}], [\ref{bib:biss2008}]).

However, just as applying the construction to reals loses the property of
ordering, more properties familiar from real and complex numbers vanish with
increasing dimension. The quaternions are only a skew field, i.e. for some $%
x,y$ we have $x$\textperiodcentered $y\neq y$\textperiodcentered $x$ for two
quaternions, the multiplication of octonions fails (in addition to not being
commutative) to be associative: for some $x,y,z$ we have $(x.y).z\neq
x.(y.z) $.

Real numbers, complex numbers, quaternions and octonions are all normed
division algebras over $\mathbb{R}$. However, by Hurwitz's theorem, they are
the only ones. The next step in the Cayley--Dickson precess, the sedenions,
in fact fails to have this structure. $\mathbb{S}$ is not a division algebra
because it has zero divisors; this means that two non-zero sedenions can be
multiplied to obtain zero: an example is\ $(e_{3}+e_{10})(e_{6}-e_{15})=0$
and the other example is $(e_{2}-e_{14})(e_{3}+e_{15})=0,$ for details see [%
\ref{bib:cawagas2004}].

In the following, we explain this doubling process.

The Cayley-Dickson algebras are a sequence $A_{0},A_{1},...$ of
non-associative $\mathbb{R}$-algebras with involution. The term
\textquotedblleft conjugation\textquotedblright\ is used to refer to the
involution because it generalizes the usual conjugation on the complex
numbers. See [\ref{bib:biss2008}] for a full explanation of the basic
properties of Cayley-Dickson algebras, We define Cayley-Dickson algebras
inductively. We begin by defining $A_{0}$ to be $\mathbb{R}.$ Given $%
A_{n-1}, $ the algebra $A_{n}$ is defined additively to be $A_{n-1}\times
A_{n-1}.$ Conjugation, multiplication and addition in $A_{n}$ are defined by%
\begin{eqnarray*}
\overline{(a,b)} &=&(\overline{a},-b), \\
(a,b)(c,d) &=&(ac-\overline{d}b,da+b\overline{c}), \\
(a,b)+(c,d) &=&(a+c,b+d),
\end{eqnarray*}%
respectively. Note that $A_{n}$ has dimension $2^{n}$ as an $\mathbb{R}-$%
vector space. If we set, as usual, $\left\Vert x\right\Vert =\sqrt{\func{Re}%
(x\overline{x})}$ for $x\in A_{n}$ then $x\overline{x}=\overline{x}%
x=\left\Vert x\right\Vert ^{2}
.$

Now, assume that $B_{16}=\{e_{i}\in \mathbb{S}:i=0,1,2,...,15\}$ is the
basis for $\mathbb{S}$, where $e_{0}$ is the identity (or unit) and $%
e_{1},e_{2},...,e_{15}$ are called imaginaries. Then a sedenion $S\in 
\mathbb{S}$ can be written as%
\begin{equation*}
S=\sum_{i=0}^{15}a_{i}e_{i}=a_{0}+\sum_{i=1}^{15}a_{i}e_{i}
\end{equation*}%
where $a_{0},a_{1},...,a_{15}$ are all real numbers. $a_{0}$ is called the
real part of $S$ and $\sum_{i=1}^{15}a_{i}e_{i}$ is called its imaginary
part.

Addition of sedenions is defined as componentwise and multiplication is
defined as follows: if $S_{1},S_{2}\in \mathbb{S}$ then we have%
\begin{equation}
S_{1}S_{2}=\left( \sum_{i=0}^{15}a_{i}e_{i}\right) \left(
\sum_{i=0}^{15}b_{i}e_{i}\right) =\sum_{i,j=0}^{15}a_{i}b_{j}(e_{i}e_{j}).
\label{equation:fvcxszeawqosbgv}
\end{equation}%
The operations requiring for the multiplication in (\ref%
{equation:fvcxszeawqosbgv}) are quite a lot. The computation of a sedenion
multiplication using the naive method requires $240$ additions and $256$
multiplications, while an algorithm which is given in [\ref{bib:cariow2013}]
can compute the same result in only $298$ additions and $122$
multiplications, see [\ref{bib:cariow2013}] for details. Moreover, efficient
algorithms for the multiplication of quaternions, octonions and
trigintaduonions with reduced number of real multiplications is already
exist in literature, see [\ref{makarov1978}], [\ref{cariow2012}] and [\ref%
{cariow2014}], respectively.

Recently, a lot of attention has been centred by theoretical physicists on
the Cayley-Dickson algebras $\mathbb{O}$ (octonions) and $\mathbb{S}$
(sedenions) because of their increasing usefulness in formulating many of
the new theories of elementary particles. In particular, the octonions $%
\mathbb{O}$ (see for example [\ref{bib:baez2002}] and [\ref{bib:okubo1995}])
has been found to be involved in so many unexpected areas (such as Clifford
algebras,\ quantum theory, topology, etc.) and sedenions appear in many
areas of science like electromagnetic theory and linear gravity.

\section{Generalized Tribonacci Sedenions and It's Generating Function and
Binet Formula}

In this section, we define generalized Tribonacci sedenions and give
generating functions and Binet formulas for them. First, we give some
information about quaternion sequences, octonion sequences and sedenion
sequences from literature.

Horadam [\ref{bib:horadam1963aa}] introduced $n$th Fibonacci and $n$th Lucas
quaternions as%
\begin{equation*}
Q_{n}=F_{n}+F_{n+1}e_{1}+F_{n+2}e_{2}+F_{n+3}e_{3}=\sum_{s=0}^{3}F_{n+s}e_{s}
\end{equation*}%
and%
\begin{equation*}
R_{n}=L_{n}+L_{n+1}e_{1}+L_{n+2}e_{2}+L_{n+3}e_{3}=\sum_{s=0}^{3}L_{n+s}e_{s}
\end{equation*}%
respectively, where $F_{n}$ and $L_{n}$ are the $n$th Fibonacci and Lucas
numbers respectively. He also defined generalized Fibonacci quaternion as%
\begin{equation*}
P_{n}=H_{n}+H_{n+1}e_{1}+H_{n+2}e_{2}+H_{n+3}e_{3}=\sum_{s=0}^{3}H_{n+s}e_{s}
\end{equation*}%
where $H_{n}$ is the $n$th generalized Fibonacci number (which is now called
Horadam number) by the recursive relation $H_{1}=p,$ $H_{2}=p+q,$ $%
H_{n}=H_{n-1}+H_{n-2}$ ($p$ and $q$ are arbitrary integers). Many other
generalization of Fibonacci quaternions has been given, see for example
Halici and Karata\c{s} [\ref{bib:halici2017}], and Polatl\i\ [\ref%
{bib:polatli2016}]. See also Tasci [\ref{tascijacobquater2017}] for
k-Jacobsthal and k-Jacobsthal-Lucas Quaternions.

Cerda-Morales [\ref{bib:cerdamorale2017}] defined and studied the
generalized Tribonacci quaternion sequence that includes the previously
introduced Tribonacci, Padovan, Narayana and third order Jacobsthal
quaternion sequences. She defined generalized Tribonacci quaternion as%
\begin{equation*}
Q_{v,n}=V_{n}+V_{n+1}e_{1}+V_{n+2}e_{2}+V_{n+3}e_{3}=%
\sum_{s=0}^{3}V_{n+s}e_{s}
\end{equation*}%
where $V_{n}$ is the $n$th generalized Tribonacci number defined by the
third-order recurrance relations%
\begin{equation*}
V_{n}=rV_{n-1}+sV_{n-2}+tV_{n-3},\text{ \ \ \ \ }
\end{equation*}%
here $V_{0}=a,V_{1}=b,V_{2}=c$ are arbitrary integers and $r,s,t$ are real
numbers. See also [\ref{akkus2017}], [\ref{szynal2017}] for Tribonacci
quaternions, [\ref{cerdamoralejacobs2017}] for third order Jacobsthal
quaternions, [\ref{cerdamoralejdualacobs2018}] for dual third order
Jacobsthal quaternions and [\ref{tasciquater2018}] for Padovan and
Pell-Padovan quaternions.

Various families of octonion number sequences (such as Fibonacci octonion,
Pell octonion, Jacobsthal octonion; and third order Jacobsthal octonion)
have been defined and studied by a number of authors in many different ways.
For example, Ke\c{c}ilioglu and Akku\c{s} [\ref{bib:kecilioglu}] introduced
the Fibonacci and Lucas octonions as%
\begin{equation*}
Q_{n}=\sum_{s=0}^{7}F_{n+s}e_{s}
\end{equation*}%
and%
\begin{equation*}
R_{n}=\sum_{s=0}^{7}L_{n+s}e_{s}
\end{equation*}%
respectively, where $F_{n}$ and $L_{n}$ are the $n$th Fibonacci and Lucas
numbers respectively. In [\ref{bib:cimen2017b}], \c{C}imen and \.{I}pek
introduced Jacobsthal octonions and Jacobsthal-Lucas octonions. In [\ref%
{bib:cerdamorale2017b}], Cerda-Morales introduced third order Jacobsthal
octonions and also in [\ref{bib:cerdamorale2018b}], she\ defined and studied
tribonacci-type octonions. See also \.{I}pek and \c{C}imen [\ref{ipekfib2016}%
] and Karata\c{s} and Halici [\ref{karatasocton2017}] for Horadam type
octonions. Furthermore, we refer to Szynal-Liana and Wlock [\ref%
{szynallianapell2016}] for Pell quaternions and Pell octonions and Catarino [%
\ref{catarinopell2016}] for Modified Pell and Modified k-Pell Quaternions
and Octonions.

A number of authors have been defined and studied sedenion number sequences
(such as second order sedenions: Fibonacci sedenion, k-Pell and
k-Pell--Lucas sedenions, Jacobsthal and Jacobsthal-Lucas sedenions). For
example, Bilgici, Toke\c{s}er and \"{U}nal [\ref{bib:bilgici2017}]
introduced the Fibonacci and Lucas sedenions as%
\begin{equation*}
\widehat{F}_{n}=\sum_{s=0}^{15}F_{n+s}e_{s}
\end{equation*}%
and%
\begin{equation*}
\widehat{L}_{n}=\sum_{s=0}^{15}L_{n+s}e_{s}
\end{equation*}%
respectively, where $F_{n}$ and $L_{n}$ are the $n$th Fibonacci and Lucas
numbers respectively. In [\ref{bib:catarino2018}], Catarino introduced
k-Pell and k-Pell--Lucas sedenions. In [\ref{bib:cimen2017}], \c{C}imen and 
\.{I}pek introduced Jacobsthal and Jacobsthal-Lucas sedenions. In [\ref%
{bib:soykansedenion2019}], Soykan introduced the $n$th Tribonacci sedenion as%
\begin{equation*}
\widehat{T}_{n}=\sum_{s=0}^{15}T_{n+s}e_{s}=T_{n}+%
\sum_{s=1}^{15}T_{n+s}e_{s}.
\end{equation*}%
(see also preprint version [\ref{soykantrisedenion2018}] of it).

G\"{u}l [\ref{bib:gul2018}] introduced the k-Fibonacci and k-Lucas
trigintaduonions as%
\begin{equation*}
TF_{k,n}=\sum_{s=0}^{31}F_{k,n+s}e_{s}
\end{equation*}%
and%
\begin{equation*}
TL_{k,n}=\sum_{s=0}^{31}L_{k,n+s}e_{s}
\end{equation*}%
respectively, where $F_{k,n}$ and $L_{k,n}$ are the $n$th k-Fibonacci and
k-Lucas numbers respectively.

We now define the generalized Tribonacci sedenion over the sedenion algebra $%
\mathbb{S}$. The $n$th generalized Tribonacci sedenion is%
\begin{equation}
\widehat{V}_{n}=\sum_{s=0}^{15}V_{n+s}e_{s}=V_{n}+\sum_{s=1}^{15}V_{n+s}e_{s}
\label{equation:bnuyvbcxdszaerm}
\end{equation}%
where $V_{n}$ is the $n$th generalized Tribonacci number defined by Eq.(\ref%
{equation:nbvcftyopuybnaeuo}).

For all non negative integer $n,$ it can be easily shown that the following
iterative relation holds: 
\begin{equation}
\widehat{V}_{n}=r\widehat{V}_{n-1}+s\widehat{V}_{n-2}+t\widehat{V}_{n-3}.
\label{equation:mnouaedcxzqwoup}
\end{equation}

The sequence $\{\widehat{V}_{n}\}_{n\geq 0}$ can be defined for negative
values of $n$ by using the recurrence (\ref{equation:mnouaedcxzqwoup}) to
extend the sequence backwards, or equivalently, by using the recurrence 
\begin{equation*}
\widehat{V}_{-n}=-\frac{s}{t}\widehat{V}_{-(n-1)}-\frac{r}{t}\widehat{V}%
_{-(n-2)}+\frac{1}{t}\widehat{V}_{-(n-3)}.
\end{equation*}%
Thus, the recurrence (\ref{equation:mnouaedcxzqwoup}) holds for all integer $%
n.$

The conjugate of $\widehat{V}_{n}$ is defined by%
\begin{equation*}
\overline{\widehat{V}_{n}}=V_{n}-%
\sum_{s=1}^{15}V_{n+s}e_{s}=V_{n}-V_{n+1}e_{1}-V_{n+2}e_{2}-...-V_{n+15}e_{15}.
\end{equation*}%
The norm of $n$th generalized Tribonacci\ sedenion is 
\begin{equation*}
\left\Vert \widehat{V}_{n}\right\Vert ^{2}=N^{2}(\widehat{V}_{n})=\widehat{V}%
_{n}\overline{\widehat{V}_{n}}=\overline{\widehat{V}_{n}}\widehat{V}%
_{n}=V_{n}^{2}+V_{n+1}^{2}+...+V_{n+15}^{2}.
\end{equation*}

Now, we will state Binet's formula for the generalized Tribonacci sedenion.
From now on, we fixed the following notations.%
\begin{equation*}
\widehat{\alpha }=\sum_{s=0}^{15}\alpha ^{s}e_{s},\text{ \ }\widehat{\beta }%
=\sum_{s=0}^{15}\beta ^{s}e_{s},\text{ \ }\widehat{\gamma }%
=\sum_{s=0}^{15}\gamma ^{s}e_{s}.
\end{equation*}

\begin{theorem}
\label{theorem:yunhbvgcfosea}For any non-negative integer $n,$ the $n$th
generalized Tribonacci sedenion is%
\begin{equation}
\widehat{V}_{n}=\frac{P\widehat{\alpha }\alpha ^{n}}{(\alpha -\beta )(\alpha
-\gamma )}+\frac{Q\widehat{\beta }\beta ^{n}}{(\beta -\alpha )(\beta -\gamma
)}+\frac{R\widehat{\gamma }\gamma ^{n}}{(\gamma -\alpha )(\gamma -\beta )}
\label{equation:fgvbvcxsdaewqzxd}
\end{equation}%
where $P,Q$ and $R$ are as in Eq.(\ref{equat:mnopcvbedcxzsa}).
\end{theorem}

\textit{Proof.} Repeated use of (\ref{equat:mnopcvbedcxzsa}) in (\ref%
{equation:bnuyvbcxdszaerm}) enable us to write$:$%
\begin{eqnarray*}
\widehat{V}_{n} &=&\sum_{s=0}^{15}V_{n+s}e_{s}=\frac{P\alpha ^{n}}{(\alpha
-\beta )(\alpha -\gamma )}+\frac{Q\beta ^{n}}{(\beta -\alpha )(\beta -\gamma
)}+\frac{R\gamma ^{n}}{(\gamma -\alpha )(\gamma -\beta )} \\
&&+\left( \frac{P\alpha ^{n+1}}{(\alpha -\beta )(\alpha -\gamma )}+\frac{%
Q\beta ^{n+1}}{(\beta -\alpha )(\beta -\gamma )}+\frac{R\gamma ^{n+1}}{%
(\gamma -\alpha )(\gamma -\beta )}\right) e_{1} \\
&&+\left( \frac{P\alpha ^{n+2}}{(\alpha -\beta )(\alpha -\gamma )}+\frac{%
Q\beta ^{n+2}}{(\beta -\alpha )(\beta -\gamma )}+\frac{R\gamma ^{n+2}}{%
(\gamma -\alpha )(\gamma -\beta )}\right) e_{2} \\
&&\vdots  \\
&&+\left( \frac{P\alpha ^{n+15}}{(\alpha -\beta )(\alpha -\gamma )}+\frac{%
Q\beta ^{n+15}}{(\beta -\alpha )(\beta -\gamma )}+\frac{R\gamma ^{n+15}}{%
(\gamma -\alpha )(\gamma -\beta )}\right) e_{15} \\
&=&\frac{P\widehat{\alpha }\alpha ^{n}}{(\alpha -\beta )(\alpha -\gamma )}+%
\frac{Q\widehat{\beta }\beta ^{n}}{(\beta -\alpha )(\beta -\gamma )}+\frac{R%
\widehat{\gamma }\gamma ^{n}}{(\gamma -\alpha )(\gamma -\beta )}.
\end{eqnarray*}%
\endproof%

We can give Binet's formula of the generalized Tribonacci sedenion for the
negative subscripts:\ for $n=1,2,3,...$ we have 
\begin{equation*}
\widehat{V}_{-n}=\frac{\allowbreak \alpha ^{2}-r\alpha -s}{t}\frac{P\widehat{%
\alpha }\alpha ^{1-n}}{(\alpha -\beta )(\alpha -\gamma )}+\frac{\allowbreak
\beta ^{2}-r\beta -s}{t}\frac{Q\widehat{\beta }\beta ^{1-n}}{(\beta -\alpha
)(\beta -\gamma )}+\frac{\allowbreak \gamma ^{2}-r\gamma -s}{t}\frac{R%
\widehat{\gamma }\gamma ^{1-n}}{(\gamma -\alpha )(\gamma -\beta )}.
\end{equation*}%
The next theorem gives us an alternatif proof of the Binet's formula for the
generalized Tribonacci sedenion. For this, we need the quadratik
approximation of $\{V_{n}\}_{n\geq 0}$:

\begin{lemma}
\lbrack \ref{bib:cerdamorale2018b}]For all integer $n,$ we have%
\begin{eqnarray*}
P\alpha ^{n+2} &=&\alpha ^{2}V_{n+2}+\alpha (sV_{n+1}+tV_{n})+tV_{n+1}, \\
Q\beta ^{n+2} &=&\beta ^{2}V_{n+2}+\beta (sV_{n+1}+tV_{n})+tV_{n+1}, \\
R\gamma ^{n+2} &=&\gamma ^{2}V_{n+2}+\gamma (sV_{n+1}+tV_{n})+tV_{n+1},
\end{eqnarray*}%
where $P,Q$ and $R$ are as in Eq.(\ref{equat:mnopcvbedcxzsa}).
\end{lemma}

\underline{Alternative Proof of Theorem \ref{theorem:yunhbvgcfosea}:}

Note that%
\begin{eqnarray*}
&&\alpha ^{2}\widehat{V}_{n+2}+\alpha (s\widehat{V}_{n+1}+t\widehat{V}_{n})+t%
\widehat{V}_{n+1} \\
&=&\alpha ^{2}(V_{n+2}+V_{n+3}e_{1}+...+V_{n+17}e_{15}) \\
&&+\alpha
((sV_{n+1}+tV_{n})+(sV_{n+2}+tV_{n+1})e_{1}+...+(sV_{n+16}+tV_{n+15})e_{15})
\\
&&+t(V_{n+1}+V_{n+2}e_{1}+...+V_{n+16}e_{15}) \\
&=&\alpha ^{2}V_{n+2}+\alpha (sV_{n+1}+tV_{n})+tV_{n+1}+(\alpha
^{2}V_{n+3}+\alpha (sV_{n+2}+tV_{n+1})+tV_{n+2})e_{1} \\
&&+(\alpha ^{2}V_{n+4}+\alpha (sV_{n+3}+tV_{n+2})+tV_{n+3})e_{2} \\
&&\vdots  \\
&&+(\alpha ^{2}V_{n+17}+\alpha (sV_{n+16}+tV_{n+15})+tV_{n+16})e_{15}.
\end{eqnarray*}%
From the identity $P\alpha ^{n+2}=\alpha ^{2}V_{n+2}+\alpha
(sV_{n+1}+tV_{n})+tV_{n+1}$ for $n$-th\ generalized Tribonacci number $V_{n},
$ we have%
\begin{equation}
\alpha ^{2}\widehat{V}_{n+2}+\alpha (s\widehat{V}_{n+1}+t\widehat{V}_{n})+t%
\widehat{V}_{n+1}=P\widehat{\alpha }\alpha ^{n+2}.
\label{equation:abctyunmbvcfoes}
\end{equation}%
Similarly, we obtain%
\begin{equation}
\beta ^{2}\widehat{V}_{n+2}+\beta (s\widehat{V}_{n+1}+t\widehat{V}_{n})+t%
\widehat{V}_{n+1}=Q\widehat{\beta }\beta ^{n+2}
\label{equation:byummnbvoepsdxcz}
\end{equation}%
and%
\begin{equation}
\gamma ^{2}\widehat{V}_{n+2}+\gamma (s\widehat{V}_{n+1}+t\widehat{V}_{n})+t%
\widehat{V}_{n+1}=R\widehat{\gamma }\gamma ^{n+2}.
\label{equation:cuytnbhvgfopsexza}
\end{equation}%
Substracting (\ref{equation:byummnbvoepsdxcz}) from (\ref%
{equation:abctyunmbvcfoes}), we have%
\begin{equation}
(\alpha +\beta )\widehat{V}_{n+2}+(s\widehat{V}_{n+1}+t\widehat{V}_{n})=%
\frac{P\widehat{\alpha }\alpha ^{n+2}-Q\widehat{\beta }\beta ^{n+2}}{\alpha
-\beta }.  \label{equation:somnbyusxczaeo}
\end{equation}%
Similarly, substracting (\ref{equation:cuytnbhvgfopsexza}) from (\ref%
{equation:abctyunmbvcfoes}), we obtain%
\begin{equation}
(\alpha +\gamma )\widehat{V}_{n+2}+(s\widehat{V}_{n+1}+t\widehat{V}_{n})=%
\frac{P\widehat{\alpha }\alpha ^{n+2}-R\widehat{\gamma }\gamma ^{n+2}}{%
\alpha -\gamma }.  \label{equation:gtybnvcxdsxzc}
\end{equation}%
Finally, substracting (\ref{equation:gtybnvcxdsxzc}) from (\ref%
{equation:somnbyusxczaeo}), we get%
\begin{eqnarray*}
\widehat{V}_{n+2} &=&\frac{1}{\beta -\gamma }\left( \frac{P\widehat{\alpha }%
\alpha ^{n+2}-Q\widehat{\beta }\beta ^{n+2}}{\alpha -\beta }-\frac{P\widehat{%
\alpha }\alpha ^{n+2}-R\widehat{\gamma }\gamma ^{n+2}}{\alpha -\gamma }%
\right)  \\
&=&\frac{P\widehat{\alpha }\alpha ^{n+2}}{(\alpha -\beta )(\alpha -\gamma )}-%
\frac{Q\widehat{\beta }\beta ^{n+2}}{(\alpha -\beta )(\beta -\gamma )}+\frac{%
R\widehat{\gamma }\gamma ^{n+2}}{(\gamma -\alpha )(\gamma -\beta )} \\
&=&\frac{P\widehat{\alpha }\alpha ^{n+2}}{(\alpha -\beta )(\alpha -\gamma )}+%
\frac{Q\widehat{\beta }\beta ^{n+2}}{(\beta -\alpha )(\beta -\gamma )}+\frac{%
R\widehat{\gamma }\gamma ^{n+2}}{(\gamma -\alpha )(\gamma -\beta )}.
\end{eqnarray*}%
So this proves (\ref{equation:fgvbvcxsdaewqzxd}). 
\endproof%

Next, we present generating function of $\widehat{V}_{n}$.

\begin{theorem}
The generating functions for the generalized Tribonacci sedenion is%
\begin{equation}
g(x)=\sum_{n=0}^{\infty }\widehat{V}_{n}x^{n}=\frac{\widehat{V}_{0}+(%
\widehat{V}_{1}-r\widehat{V}_{0})x+t\widehat{V}_{-1}x^{2}}{1-rx-sx^{2}-tx^{3}%
}.  \label{equation:hnbvyuedsxzapo}
\end{equation}
\end{theorem}

\textit{Proof.} Define $g(x)=\sum_{n=0}^{\infty }\widehat{V}_{n}x^{n}.$ Note
that%
\begin{equation*}
\begin{array}{ccccccccc}
g(x) & = & \widehat{V}_{0}+ & \widehat{V}_{1}x+ & \widehat{V}_{2}x^{2}+ & 
\widehat{V}_{3}x^{3}+ & \widehat{V}_{4}x^{4}+ & \widehat{V}_{5}x^{5}+...+ & 
\widehat{V}_{n}x^{n}\text{ }+... \\ 
rxg(x) & = &  & r\widehat{V}_{0}x+ & r\widehat{V}_{1}x^{2}+ & r\widehat{V}%
_{2}x^{3}+ & r\widehat{V}_{3}x^{4}+ & r\widehat{V}_{4}x^{5}+...+ & r\widehat{%
V}_{n-1}x^{n}+... \\ 
sx^{2}g(x) & = &  &  & s\widehat{V}_{0}x^{2}+ & s\widehat{V}_{1}x^{3}+ & s%
\widehat{V}_{2}x^{4}+ & s\widehat{V}_{3}x^{5}+...+ & s\widehat{V}%
_{n-2}x^{n}+... \\ 
tx^{3}g(x) & = &  &  &  & t\widehat{V}_{0}x^{3}+ & t\widehat{V}_{1}x^{4}+ & t%
\widehat{V}_{2}x^{5}+...+ & t\widehat{V}_{n-3}x^{n}+...%
\end{array}%
\end{equation*}%
Using above table and the recurans $\widehat{V}_{n}=r\widehat{V}_{n-1}+s%
\widehat{V}_{n-2}+t\widehat{V}_{n-3}$ we have%
\begin{eqnarray*}
(1-rx-sx^{2}-tx^{3})g(x) &=&g(x)-rxg(x)-sx^{2}g(x)-tx^{3}g(x) \\
&=&\widehat{V}_{0}+(\widehat{V}_{1}-r\widehat{V}_{0})x+(\widehat{V}_{2}-r%
\widehat{V}_{1}-s\widehat{V}_{0})x^{2}+(\widehat{V}_{3}-r\widehat{V}_{2}-s%
\widehat{V}_{1}-t\widehat{V}_{0})x^{3}+ \\
&&(\widehat{V}_{4}-r\widehat{V}_{3}-s\widehat{V}_{2}-t\widehat{V}%
_{1})x^{4}+...+(\widehat{V}_{n}-r\widehat{V}_{n-1}-s\widehat{V}_{n-2}-t%
\widehat{V}_{n-3}+)x^{n}+... \\
&=&\widehat{V}_{0}+(\widehat{V}_{1}-r\widehat{V}_{0})x+(\widehat{V}_{2}-r%
\widehat{V}_{1}-s\widehat{V}_{0})x^{2}.
\end{eqnarray*}%
It follows that%
\begin{equation*}
g(x)=\frac{\widehat{V}_{0}+(\widehat{V}_{1}-r\widehat{V}_{0})x+(\widehat{V}%
_{2}-r\widehat{V}_{1}-s\widehat{V}_{0})x^{2}}{1-rx-sx^{2}-tx^{3}}.
\end{equation*}%
Since $\widehat{V}_{2}-r\widehat{V}_{1}-s\widehat{V}_{0}=t\widehat{V}_{-1},$
the generating function for the Tribonacci sedenion is%
\begin{equation*}
g(x)=\frac{\widehat{V}_{0}+(\widehat{V}_{1}-r\widehat{V}_{0})x+t\widehat{V}%
_{-1}x^{2}}{1-rx-sx^{2}-tx^{3}}.
\end{equation*}%
\endproof%

In the following theorem we present another form of Binet formula for the
generalized Tribonacci sedenion using generating function.

\begin{theorem}
\label{theorem:vbghcxdszrtupomn}For any integer $n,$ the $n$th generalized
Tribonacci sedenion is%
\begin{eqnarray*}
\widehat{V}_{n} &=&\frac{((\alpha ^{2}-r\alpha )\widehat{V}_{0}+\alpha 
\widehat{V}_{1}+t\widehat{V}_{-1})\alpha ^{n}}{\left( \alpha -\gamma \right)
\left( \alpha -\beta \right) }+\frac{((\beta ^{2}-r\beta )\widehat{V}%
_{0}+\beta \widehat{V}_{1}+t\widehat{V}_{-1})\beta ^{n}}{\left( \beta
-\gamma \right) \left( \beta -\alpha \right) } \\
&&+\frac{((\gamma ^{2}-r\gamma )\widehat{V}_{0}+\gamma \widehat{V}_{1}+t%
\widehat{V}_{-1})\gamma ^{n}}{\left( \gamma -\alpha \right) \left( \gamma
-\beta \right) }.
\end{eqnarray*}
\end{theorem}

\textit{Proof.} We can use generating function. Since%
\begin{equation*}
1-rx-sx^{2}-tx^{3}=(1-\alpha x)(1-\beta x)(1-\gamma x)
\end{equation*}%
we can write the generating function of $\widehat{V}_{n}$ as \ 
\begin{eqnarray*}
g(x) &=&\frac{\widehat{V}_{0}+(\widehat{V}_{1}-r\widehat{V}_{0})x+t\widehat{V%
}_{-1}x^{2}}{1-rx-sx^{2}-tx^{3}}=\frac{\widehat{V}_{0}+(\widehat{V}_{1}-r%
\widehat{V}_{0})x+t\widehat{V}_{-1}x^{2}}{(1-\alpha x)(1-\beta x)(1-\gamma x)%
} \\
&=&\frac{D}{(1-\alpha x)}+\frac{E}{(1-\beta x)}+\frac{F}{(1-\gamma x)} \\
&=&\frac{D(1-\beta x)(1-\gamma x)+E(1-\alpha x)(1-\gamma x)+F(1-\alpha
x)(1-\beta x)}{(1-\alpha x)(1-\beta x)(1-\gamma x)} \\
&=&\frac{(D+E+F)+(-D\beta -D\gamma -E\alpha -E\gamma -F\alpha -F\beta
)x+(D\beta \gamma +E\alpha \gamma +F\alpha \beta )x^{2}}{(1-\alpha
x)(1-\beta x)(1-\gamma x)}
\end{eqnarray*}%
We need to find $D,E$ and $F$, so the following system of equations should
be solved:%
\begin{eqnarray*}
D+E+F &=&\widehat{V}_{0} \\
-D\beta -D\gamma -E\alpha -E\gamma -F\alpha -F\beta &=&\widehat{V}_{1}-r%
\widehat{V}_{0} \\
D\beta \gamma +E\alpha \gamma +F\alpha \beta &=&t\widehat{V}_{-1}
\end{eqnarray*}%
We find that%
\begin{eqnarray*}
D &=&\frac{t\widehat{V}_{-1}+\alpha \widehat{V}_{1}+\alpha ^{2}\widehat{V}%
_{0}-r\alpha \widehat{V}_{0}}{\alpha ^{2}-\alpha \beta -\alpha \gamma +\beta
\gamma }=\frac{(\alpha ^{2}-r\alpha )\widehat{V}_{0}+\alpha \widehat{V}_{1}+t%
\widehat{V}_{-1}}{\left( \alpha -\gamma \right) \left( \alpha -\beta \right) 
} \\
E &=&\frac{t\widehat{V}_{-1}+\beta \widehat{V}_{1}+\beta ^{2}\widehat{V}%
_{0}-r\beta V_{0}}{\beta ^{2}-\alpha \beta +\alpha \gamma -\beta \gamma }=%
\frac{(\beta ^{2}-r\beta )\widehat{V}_{0}+\beta \widehat{V}_{1}+t\widehat{V}%
_{-1}}{\left( \beta -\gamma \right) \left( \beta -\alpha \right) } \\
F &=&\frac{t\widehat{V}_{-1}+\gamma \widehat{V}_{1}+\gamma ^{2}\widehat{V}%
_{0}-r\gamma \widehat{V}_{0}}{\gamma ^{2}+\alpha \beta -\alpha \gamma -\beta
\gamma }=\frac{(\gamma ^{2}-r\gamma )\widehat{V}_{0}+\gamma \widehat{V}_{1}+t%
\widehat{V}_{-1}}{\left( \gamma -\alpha \right) \left( \gamma -\beta \right) 
}
\end{eqnarray*}%
and 
\begin{eqnarray*}
g(x) &=&\frac{(\alpha ^{2}-r\alpha )\widehat{V}_{0}+\alpha \widehat{V}_{1}+t%
\widehat{V}_{-1}}{\left( \alpha -\gamma \right) \left( \alpha -\beta \right) 
}\sum_{n=0}^{\infty }\alpha ^{n}x^{n}+\frac{(\beta ^{2}-r\beta )\widehat{V}%
_{0}+\beta \widehat{V}_{1}+t\widehat{V}_{-1}}{\left( \beta -\gamma \right)
\left( \beta -\alpha \right) }\sum_{n=0}^{\infty }\beta ^{n}x^{n} \\
&&+\frac{(\gamma ^{2}-r\gamma )\widehat{V}_{0}+\gamma \widehat{V}_{1}+t%
\widehat{V}_{-1}}{\left( \gamma -\alpha \right) \left( \gamma -\beta \right) 
}\sum_{n=0}^{\infty }\gamma ^{n}x^{n} \\
&=&\sum_{n=0}^{\infty }\left( 
\begin{array}{c}
\dfrac{((\alpha ^{2}-r\alpha )\widehat{V}_{0}+\alpha \widehat{V}_{1}+t%
\widehat{V}_{-1})\alpha ^{n}}{\left( \alpha -\gamma \right) \left( \alpha
-\beta \right) }+\dfrac{((\beta ^{2}-r\beta )\widehat{V}_{0}+\beta \widehat{V%
}_{1}+t\widehat{V}_{-1})\beta ^{n}}{\left( \beta -\gamma \right) \left(
\beta -\alpha \right) } \\ 
+\dfrac{((\gamma ^{2}-r\gamma )\widehat{V}_{0}+\gamma \widehat{V}_{1}+t%
\widehat{V}_{-1})\gamma ^{n}}{\left( \gamma -\alpha \right) \left( \gamma
-\beta \right) }%
\end{array}%
\right) x^{n}.
\end{eqnarray*}%
Thus Binet formula of Tribonacci sedenion is 
\begin{eqnarray*}
\widehat{V}_{n} &=&\frac{((\alpha ^{2}-r\alpha )\widehat{V}_{0}+\alpha 
\widehat{V}_{1}+t\widehat{V}_{-1})\alpha ^{n}}{\left( \alpha -\gamma \right)
\left( \alpha -\beta \right) }+\frac{((\beta ^{2}-r\beta )\widehat{V}%
_{0}+\beta \widehat{V}_{1}+t\widehat{V}_{-1})\beta ^{n}}{\left( \beta
-\gamma \right) \left( \beta -\alpha \right) } \\
&&+\frac{((\gamma ^{2}-r\gamma )\widehat{V}_{0}+\gamma \widehat{V}_{1}+t%
\widehat{V}_{-1})\gamma ^{n}}{\left( \gamma -\alpha \right) \left( \gamma
-\beta \right) }.
\end{eqnarray*}%
\endproof%

If we compare Theorem \ref{theorem:yunhbvgcfosea} and Theorem \ref%
{theorem:vbghcxdszrtupomn} and use the definition of $\widehat{V}_{n},$ we
have the following Remark showing relations between $\widehat{T}_{-1},%
\widehat{T}_{0},\widehat{T}_{1}\ $and $\widehat{\alpha },\widehat{\beta },%
\widehat{\gamma }.$ We obtain (b) after solving the system of the equations
in (a).

\begin{remark}
We have the following identities:

\begin{description}
\item[(a)] 
\begin{eqnarray*}
((\alpha ^{2}-r\alpha )\widehat{V}_{0}+\alpha \widehat{V}_{1}+t\widehat{V}%
_{-1}) &=&P\widehat{\alpha } \\
((\beta ^{2}-r\beta )\widehat{V}_{0}+\beta \widehat{V}_{1}+t\widehat{V}%
_{-1}) &=&Q\widehat{\beta } \\
((\gamma ^{2}-r\gamma )\widehat{V}_{0}+\gamma \widehat{V}_{1}+t\widehat{V}%
_{-1}) &=&R\widehat{\gamma }
\end{eqnarray*}

\item[(b)] 
\begin{eqnarray*}
\sum_{s=0}^{15}V_{-1+s}e_{s} &=&\widehat{V}_{-1}=\frac{P\beta \gamma (\gamma
-\beta )\widehat{\alpha }+Q\alpha \gamma (\alpha -\gamma )\widehat{\beta }%
+R\alpha \beta (\beta -\alpha )\widehat{\gamma }}{t\left( \gamma -\beta
\right) \left( \alpha -\gamma \right) \left( \alpha -\beta \right) } \\
\sum_{s=0}^{15}V_{s}e_{s} &=&\widehat{V}_{0}=\frac{P(\beta -\gamma )\widehat{%
\alpha }+Q(\gamma -\alpha )\widehat{\beta }+R(\alpha -\beta )\widehat{\gamma 
}}{\left( \beta -\gamma \right) \left( \alpha -\gamma \right) \left( \alpha
-\beta \right) } \\
\sum_{s=0}^{15}V_{1+s}e_{s} &=&\widehat{V}_{1}=\dfrac{\widehat{\alpha }E_{1}+%
\widehat{\beta }E_{2}+\widehat{\gamma }E_{3}}{\left( \gamma -\beta \right)
\left( \alpha -\gamma \right) \left( \alpha -\beta \right) }
\end{eqnarray*}%
where $E_{1}=P\left( \beta -\gamma \right) \left( -r+\beta +\gamma \right)
,E_{2}=Q\left( \gamma -\alpha \right) \left( -r+\alpha +\gamma \right)
,E_{3}=R\left( \alpha -\beta \right) \left( -r+\alpha +\beta \right) .$
\end{description}
\end{remark}

Using above Remark we can find $\widehat{V}_{2}$ as follows:%
\begin{equation*}
\sum_{s=0}^{15}V_{2+s}e_{s}=\widehat{V}_{2}=r\widehat{V}_{1}+s\widehat{V}%
_{0}+t\widehat{V}_{-1}=\frac{\widehat{\alpha }c_{1}+\widehat{\beta }c_{2}+%
\widehat{\gamma }c_{3}}{\left( \beta -\gamma \right) \left( \alpha -\gamma
\right) \left( \alpha -\beta \right) }
\end{equation*}%
where 
\begin{eqnarray*}
c_{1} &=&P\left( \beta -\gamma \right) \left( s-r\beta -r\gamma +\beta
\gamma +r^{2}\right) , \\
c_{2} &=&Q\left( \gamma -\alpha \right) \left( s-r\alpha -r\gamma +\alpha
\gamma +r^{2}\right) , \\
c_{3} &=&R\left( \alpha -\beta \right) \left( s-r\alpha -r\beta +\alpha
\beta +r^{2}\right) .
\end{eqnarray*}

Now, we present the formula which gives the summation of the first $n$
generalized Tribonacci numbers.

\begin{lemma}[{[\protect\ref{bib:cerdamorale2017}]}]
For every integer $n\geq 0,$ we have 
\begin{equation}
\sum\limits_{i=0}^{n}V_{i}=\dfrac{%
V_{n+2}+(1-r)V_{n+1}+tV_{n}+(r+s-1)V_{0}+(r-1)V_{1}-V_{2}}{\varepsilon
_{r,s,t}}  \label{equation:easzpuyhbnvcdf}
\end{equation}%
where $\varepsilon =\varepsilon _{r,s,t}=r+s+t-1$.
\end{lemma}

Next, we present the formula which gives the summation of the first $n$
generalized Tribonacci sedenion.

\begin{theorem}
The summation formula for the generalized Tribonacci sedenion is%
\begin{equation}
\sum\limits_{i=0}^{n}\widehat{V}_{i}=\dfrac{\widehat{V}_{n+2}+(1-r)\widehat{V%
}_{n+1}+t\widehat{V}_{n}+\varphi _{r,s,t}}{\varepsilon _{r,s,t}}
\label{equation:yaomazxtyhnbvgfc}
\end{equation}%
where $\varepsilon =\varepsilon _{r,s,t}=r+s+t-1,$ $\mu =\mu
_{r,s,t}=(r+s-1)V_{0}+(r-1)V_{1}-V_{2}$ and $\varphi =\varphi _{r,s,t}=$ $%
\mu +e_{1}(\mu -\varepsilon V_{0})+e_{2}(\mu -\varepsilon
(V_{0}+V_{1})+....+e_{15}(\mu -\varepsilon (V_{0}+V_{1}+...+V_{14})).$
\end{theorem}

\textit{Proof.} Using (\ref{equation:bnuyvbcxdszaerm}) we obtain%
\begin{eqnarray*}
\sum\limits_{i=0}^{n}\widehat{V}_{i}
&=&\sum\limits_{i=0}^{n}V_{i}+e_{1}\sum\limits_{i=0}^{n}V_{i+1}+e_{2}\sum%
\limits_{i=0}^{n}V_{i+2}+...+e_{15}\sum\limits_{i=0}^{n}V_{i+15} \\
&=&(V_{0}+...+V_{n})+e_{1}(V_{1}+...+V_{n+1}) \\
&&+e_{2}(V_{2}+...+V_{n+2})+...+e_{15}(V_{15}+...+V_{n+15}).
\end{eqnarray*}%
Using (\ref{equation:easzpuyhbnvcdf}) and the notation $\mu
_{r,s,t}=(r+s-1)V_{0}+(r-1)V_{1}-V_{2},$ it can be written that 
\begin{eqnarray*}
\varepsilon _{r,s,t}\sum\limits_{i=0}^{n}\widehat{V}_{i}
&=&(V_{n+2}+(1-r)V_{n+1}+tV_{n}+\mu _{r,s,t}) \\
&&+e_{1}(V_{n+3}+(1-r)V_{n+2}+tV_{n+1}+\mu -\varepsilon V_{0}) \\
&&+e_{2}(V_{n+4}+(1-r)V_{n+3}+tV_{n+2}+\mu -\varepsilon (V_{0}+V_{1})) \\
&&\vdots \\
&&+e_{15}(V_{n+17}+(1-r)V_{n+16}+tV_{n+15}+\mu -\varepsilon
(V_{0}+V_{1}+...+V_{14})) \\
&=&\widehat{V}_{n+2}+(1-r)\widehat{V}_{n+1}+t\widehat{V}_{n}+\varphi _{r,s,t}
\end{eqnarray*}%
where $\varphi =\varphi _{r,s,t}=$ $\mu +e_{1}(\mu -\varepsilon
V_{0})+e_{2}(\mu -\varepsilon (V_{0}+V_{1})+....+e_{15}(\mu -\varepsilon
(V_{0}+V_{1}+...+V_{14})).$ It now follows that%
\begin{equation*}
\sum\limits_{i=0}^{n}\widehat{V}_{i}=\dfrac{\widehat{V}_{n+2}+(1-r)\widehat{V%
}_{n+1}+t\widehat{V}_{n}+\varphi _{r,s,t}}{\varepsilon _{r,s,t}}
\end{equation*}%
and this completes the proof. 
\endproof%

The formula (\ref{equation:cvrtfgouytghfds}) of next Theorem gives the norm
of the generalized Tribonacci sedenion.

\begin{theorem}
The norm of generalized Tribonacci sedenion is given by%
\begin{equation}
\left\Vert \widehat{V}_{n}\right\Vert ^{2}=\frac{(\beta -\alpha )^{2}P^{2}%
\widetilde{\alpha }\alpha ^{2n}+(\alpha -\gamma )^{2}Q^{2}\widetilde{\beta }%
\beta ^{2n}+(\alpha -\beta )^{2}R^{2}\widetilde{\gamma }\gamma ^{2n}-2M}{%
\psi ^{2}}  \label{equation:cvrtfgouytghfds}
\end{equation}%
where 
\begin{eqnarray}
\psi &=&\psi (\alpha ,\beta ,\gamma )=(\alpha -\beta )(\alpha -\gamma
)(\beta -\gamma ),  \notag \\
\widetilde{\alpha } &=&1+\alpha ^{2}+\alpha ^{4}+...+\alpha ^{30},  \notag \\
\widetilde{\beta } &=&1+\beta ^{2}+\beta ^{4}+...+\beta ^{30},  \notag \\
\widetilde{\gamma } &=&1+\gamma ^{2}+\gamma ^{4}+...+\gamma ^{30},
\label{equation:erwdcvsxfaosghbvc} \\
\widetilde{\alpha \beta } &=&1+\alpha \beta +(\alpha \beta )^{2}+(\alpha
\beta )^{3}+...+(\alpha \beta )^{15},  \notag \\
\widetilde{\alpha \gamma } &=&1+\alpha \gamma +(\alpha \gamma )^{2}+(\alpha
\gamma )^{3}+...+(\alpha \gamma )^{15},  \notag \\
\widetilde{\beta \gamma } &=&1+\beta \gamma +(\beta \gamma )^{2}+(\beta
\gamma )^{3}+...+(\beta \gamma )^{15},  \notag \\
M &=&(\alpha -\gamma )(\beta -\gamma )PQ\widetilde{\alpha \beta }(\alpha
\beta )^{n}+(\alpha -\beta )(\beta -\gamma )PQ\widetilde{\alpha \gamma }%
(\alpha \gamma )^{n}+(\alpha -\beta )(\alpha -\gamma )PQ\widetilde{\beta
\gamma }(\beta \gamma )^{n}.  \notag
\end{eqnarray}
\end{theorem}

\textit{Proof.} We know that $\left\Vert \widehat{V}_{n}\right\Vert
^{2}=\sum_{s=0}^{15}V_{n+s}^{2}.$ Using Binet formula of $V_{n}$ we have%
\begin{equation*}
\psi V_{n}=(\beta -\gamma )P\alpha ^{n}+(\gamma -\alpha )Q\beta ^{n}+(\alpha
-\beta )R\gamma ^{n}
\end{equation*}%
where $\psi $ is as in (\ref{equation:erwdcvsxfaosghbvc})$.$ It follows that%
\begin{eqnarray*}
\psi ^{2}V_{n}^{2} &=&(\beta -\gamma )^{2}P^{2}\alpha ^{2n}+(\gamma -\alpha
)^{2}Q^{2}\beta ^{2n} \\
&&+(\alpha -\beta )^{2}R^{2}\gamma ^{2n}+2(\beta -\gamma )(\gamma -\alpha
)PQ(\alpha \beta )^{n} \\
&&+2(\beta -\gamma )(\alpha -\beta )PR(\alpha \gamma )^{n}+2(\gamma -\alpha
)(\alpha -\beta )QR(\beta \gamma )^{n}
\end{eqnarray*}%
and%
\begin{eqnarray*}
\psi ^{2}\left\Vert \widehat{V}_{n}\right\Vert ^{2} &=&\psi
^{2}(V_{n}^{2}+V_{n+1}^{2}+V_{n+2}^{2}+...+V_{n+15}^{2}) \\
&=&(\beta -\gamma )^{2}P^{2}\widetilde{\alpha }\alpha ^{2n}+(\gamma -\alpha
)^{2}Q^{2}\widetilde{\beta }\beta ^{2n} \\
&&+(\alpha -\beta )^{2}R^{2}\widetilde{\gamma }\gamma ^{2n}+2(\beta -\gamma
)(\gamma -\alpha )PQ\widetilde{\alpha \beta }(\alpha \beta )^{n} \\
&&+2(\beta -\gamma )(\alpha -\beta )PR\widetilde{\alpha \gamma }(\alpha
\gamma )^{n}+2(\gamma -\alpha )(\alpha -\beta )QR\widetilde{\beta \gamma }%
(\beta \gamma )^{n}
\end{eqnarray*}%
where $\widetilde{\alpha },\widetilde{\beta },\widetilde{\gamma },\widetilde{%
\alpha \beta },\widetilde{\alpha \gamma }$ and $\widetilde{\beta \gamma }$
are as in (\ref{equation:erwdcvsxfaosghbvc}). 
\endproof%

Now, we present the quadratic approximation of $\widehat{V}_{n}$.

\begin{theorem}
For all integer $n\geq 0,$ we have%
\begin{eqnarray*}
P\widehat{\alpha }\alpha ^{n+2} &=&\alpha ^{2}\widehat{V}_{n+2}+\alpha (s%
\widehat{V}_{n+1}+t\widehat{V}_{n})+t\widehat{V}_{n+1}, \\
Q\widehat{\beta }\beta ^{n+2} &=&\beta ^{2}\widehat{V}_{n+2}+\beta (s%
\widehat{V}_{n+1}+t\widehat{V}_{n})+t\widehat{V}_{n+1}, \\
R\widehat{\gamma }\gamma ^{n+2} &=&\gamma ^{2}\widehat{V}_{n+2}+\gamma (s%
\widehat{V}_{n+1}+t\widehat{V}_{n})+t\widehat{V}_{n+1},
\end{eqnarray*}%
where $P,Q$ and $R$ are as in (\ref{equat:mnopcvbedcxzsa}).
\end{theorem}

\textit{Proof. }\ Using Binet's formula of $\widehat{V}_{n},$ we obtain 
\begin{eqnarray*}
\alpha \widehat{V}_{n+2}+(s+\beta \gamma )\widehat{V}_{n+1}+t\widehat{V}_{n}
&=&\frac{P\widehat{\alpha }\alpha ^{n}(\alpha ^{3}+(s+\beta \gamma )\alpha
+t)}{(\alpha -\beta )(\alpha -\gamma )} \\
&&+\frac{Q\widehat{\beta }\beta ^{n}(\alpha \beta ^{2}+(s+\beta \gamma
)\beta +t)}{(\beta -\alpha )(\beta -\gamma )} \\
&&+\frac{R\widehat{\gamma }\gamma ^{n}(\alpha \gamma ^{2}+(s+\beta \gamma
)\gamma +t)}{(\gamma -\alpha )(\gamma -\beta )}.
\end{eqnarray*}%
Since $\alpha \beta ^{2}+(s+\beta \gamma )\beta +t=0$, $\alpha \gamma
^{2}+(s+\beta \gamma )\gamma +t=0$ and 
\begin{equation*}
\alpha \left( \alpha -\beta \right) \left( \alpha -\gamma \right)
=\allowbreak \alpha ^{3}-\alpha ^{2}\beta -\alpha ^{2}\gamma +\alpha \beta
\gamma =\allowbreak \alpha ^{3}-(\alpha \beta +\alpha \gamma )\alpha
+t=\allowbreak \allowbreak \alpha ^{3}+(s+\beta \gamma )\alpha +t
\end{equation*}%
it follows that%
\begin{equation}
\alpha \widehat{V}_{n+2}+(s+\beta \gamma )\widehat{V}_{n+1}+t\widehat{V}%
_{n}=P\widehat{\alpha }\alpha ^{n+1}.  \label{equat:ybnvuerdxszaeomn}
\end{equation}%
Multiplying (\ref{equat:ybnvuerdxszaeomn}) by $\alpha $ and using again $%
\alpha \beta \gamma =t,$ we obtain 
\begin{equation*}
P\widehat{\alpha }\alpha ^{n+2}=\alpha ^{2}\widehat{V}_{n+2}+\alpha (s+\beta
\gamma )\widehat{V}_{n+1}+\alpha t\widehat{V}_{n}=\alpha ^{2}\widehat{V}%
_{n+2}+\alpha (s\widehat{V}_{n+1}+t\widehat{V}_{n})+t\widehat{V}_{n+1}.
\end{equation*}%
If we change the role of $\alpha $ with $\beta $ and $\gamma ,$ we obtain
the desired result. 
\endproof%

\section{Some Identities For Special Cases of The Generalized Tribonacci
Sedenion}

In this section, we give identities about Tribonacci and Tribonacci-Lucas
sedenions.

\begin{theorem}
\label{theorem:bnutyvbcdsaerbv}For $n\geq 1,$ the following identities hold:

\begin{description}
\item[(a)] $\widehat{K}_{n}=3\widehat{T}_{n+1}-2\widehat{T}_{n}-\widehat{T}%
_{n-1},$ here $T_{n}=V(0,1,1;1,1,1)$ and $K_{n}=V(3,1,3;1,1,1)$.

\item[(b)] $\widehat{V}_{n}+\overline{\widehat{V}_{n}}=2V_{n}$

\item[(c)] $\widehat{V}_{n+1}+\widehat{V}_{n}=\dfrac{P\widehat{\alpha }%
\left( \alpha +1\right) \alpha ^{n}}{(\alpha -\beta )(\alpha -\gamma )}+%
\dfrac{Q\widehat{\beta }\left( \beta +1\right) \beta ^{n}}{(\beta -\alpha
)(\beta -\gamma )}+\dfrac{R\widehat{\gamma }\left( \gamma +1\right) \gamma
^{n}}{(\gamma -\alpha )(\gamma -\beta )},$

\item[(d)] $\sum\limits_{i=0}^{n}\binom{n}{i}\widehat{V}_{i}=\dfrac{P%
\widehat{\alpha }\alpha }{(\alpha -\beta )(\alpha -\gamma )}(1+\alpha )^{n}+%
\dfrac{Q\widehat{\beta }\beta }{(\beta -\alpha )(\beta -\gamma )}(1+\beta
)^{n}+\dfrac{R\widehat{\gamma }\gamma }{(\gamma -\alpha )(\gamma -\beta )}%
(1+\gamma )^{n},$
\end{description}
\end{theorem}

\textit{Proof. }Since $K_{n}=3T_{n+1}-2T_{n}-T_{n-1}$, (a) follows. The
others can be established easily.\ \ 
\endproof%

\begin{theorem}
For $n\geq 0,$ $m\geq 3$, we have

\begin{description}
\item[(a)] $\widehat{K}_{m+n}=K_{n-1}\widehat{T}_{m+2}+(\widehat{T}_{m+1}+%
\widehat{T}_{m})K_{n-2}+K_{n-3}\widehat{T}_{m+1},$

\item[(b)] $\widehat{K}_{m+n}=K_{m+2}\widehat{T}_{n-1}+(K_{m+1}+K_{m})%
\widehat{T}_{n-2}+K_{m+1}\widehat{T}_{n-3}.$
\end{description}
\end{theorem}

\textit{Proof.} (a) can be proved by strong induction on $n$ and (b) can be
proved by strong induction on $m$. 
\endproof%

Similar formula can be given for the other Fibonacci sequences (for example,
for $\widehat{P}_{m+n}$ and $\widehat{R}_{m+n}$).

\begin{theorem}[{[\protect\ref{bib:howard2001}]}]
Consider the sequences $J_{n}$ which is called generalized Tribonacci-Lucas
sequence, and $V_{n}$ given by 
\begin{equation*}
J_{n}=\{V_{n}(3,r,r^{2}+2s;r,s,t)\}.
\end{equation*}%
and $V_{n}=\{V_{n}(V_{0},V_{1},V_{2};r,s,t)\}.$ Then, for all integers $n$
and $m$ we have%
\begin{equation*}
V_{n+2m}=J_{m}V_{n+m}-t^{m}J_{-m}V_{n}+t^{m}V_{n-m}.
\end{equation*}
\end{theorem}

Using above theorem, we have the following corollary.

\begin{corollary}
For all integers $n$ and $m,$ the followings are true:

\begin{description}
\item[(a)] $T_{n+2m}=K_{m}T_{n+m}-K_{-m}T_{n}+T_{n-m}$

\item[(b)] $K_{n+2m}=K_{m}K_{n+m}-K_{-m}K_{n}+K_{n-m}.$
\end{description}
\end{corollary}

\textit{Proof.}

\begin{description}
\item[(a)] Take $V_{n}=\{V_{n}(0,1,1;1,1,1)\}=\{T_{n}\}$ and $%
J_{n}=\{V_{n}(3,1,3;1,1,1)\}=\{K_{n}\}.$

\item[(b)] Take $V_{n}=\{V_{n}(3,1,3;1,1,1)\}=\{K_{n}\}$ and $%
J_{n}=\{V_{n}(3,1,3;1,1,1)\}=\{K_{n}\}.$
\end{description}

\section{Matrices and Determinants Related with Tribonacci and
Tribonacci-Lucas Sedenions}

Let $\{X_{n}\}$\ and $\{Y_{n}\}$\ be any sequences generated\ by (\ref%
{equation:nbvcftyopuybnaeuo}). Define the $4\times 4$ determinant $D_{n},$
for all integers $n,$ by%
\begin{equation*}
D_{n}=\left\vert 
\begin{array}{cccc}
X_{n} & Y_{n} & Y_{n+1} & Y_{n+2} \\ 
X_{2} & Y_{2} & Y_{3} & Y_{4} \\ 
X_{1} & Y_{1} & Y_{2} & Y_{3} \\ 
X_{0} & Y_{0} & Y_{1} & Y_{2}%
\end{array}%
\right\vert .
\end{equation*}

\begin{theorem}
\label{theorem:fgstyaremncvbx}$D_{n}=0$ for all integers $n.$
\end{theorem}

\textit{Proof.} This is a special case of a result in [\ref{bib:melham1995}].

\begin{theorem}
The followings are true.

\begin{description}
\item[(a)] $44\widehat{T}_{n}=10\widehat{K}_{n+2}-6\widehat{K}_{n+1}-8%
\widehat{K}_{n}.$

\item[(b)] $\widehat{K}_{n}=-\widehat{T}_{n+2}+4\widehat{T}_{n+1}-\widehat{T}%
_{n}.$
\end{description}
\end{theorem}

\textit{Proof.} Taking $\{X_{n}\}$\ $=\{T_{n}\},\{Y_{n}\}$\ $=\{K_{n}\}$ and
applying Theorem \ref{theorem:fgstyaremncvbx} and expanding $D_{n}$ along
the top row gives $44T_{n}=10K_{n+2}-6K_{n+1}-8K_{n}\ $and now (a) follows.
Taking $\{X_{n}\}$\ $=\{K_{n}\},\{Y_{n}\}$\ $=\{T_{n}\}$ and applying
Theorem \ref{theorem:fgstyaremncvbx} and expanding $D_{n}$ along the top row
gives $K_{n}=-T_{n+2}+4T_{n+1}-T_{n}\ $and now (b) follows. 
\endproof%

Note that similar formulas can be given for the other Tribonacci sequences
(for example, taking $\{X_{n}\}$\ $=\{P_{n}\},\{Y_{n}\}$\ $=\{R_{n}\}$ and
vice versa).

Consider the sequence $\{U_{n}\}$ which is defined by the third-order
recurrence relation%
\begin{equation*}
U_{n}=U_{n-1}+U_{n-2}+U_{n-3},\text{ \ \ \ \ }U_{0}=U_{1}=0,U_{2}=1.
\end{equation*}%
Note that some authors call $\{U_{n}\}$ as a Tribonacci sequence instead of $%
\{T_{n}\}$ . The numbers $U_{n}$ can be expressed using Binet's formula%
\begin{equation*}
U_{n}=\frac{\alpha ^{n}}{(\alpha -\beta )(\alpha -\gamma )}+\frac{\beta ^{n}%
}{(\beta -\alpha )(\beta -\gamma )}+\frac{\gamma ^{n}}{(\gamma -\alpha
)(\gamma -\beta )}
\end{equation*}%
and the negative numbers $U_{-n}$ $(n=1,2,3,...)$ satisfies the recurrence
relation%
\begin{equation*}
U_{-n}=\left\vert 
\begin{array}{cc}
U_{n+1} & U_{n+2} \\ 
U_{n} & U_{n+1}%
\end{array}%
\right\vert =U_{n+1}^{2}-U_{n+2}U_{n}.
\end{equation*}%
The matrix method is very useful method in order to obtain some identities
for special sequences. We define the square matrix $M$ of order $3$ as:%
\begin{equation*}
M=\left( 
\begin{array}{ccc}
r & s & t \\ 
1 & 0 & 0 \\ 
0 & 1 & 0%
\end{array}%
\right)
\end{equation*}%
such that $\det M=t.$ Note that 
\begin{equation}
M^{n}=\left( 
\begin{array}{ccc}
U_{n+2} & sU_{n+1}+tU_{n} & tU_{n+1} \\ 
U_{n+1} & sU_{n}+tU_{n-1} & tU_{n} \\ 
U_{n} & sU_{n-1}+tU_{n-2} & tU_{n-1}%
\end{array}%
\right) .  \label{equation:gbhasmnuopdcxsz}
\end{equation}%
For a proof of (\ref{equation:gbhasmnuopdcxsz}), see [\ref{bib:basu:2014}].
Matrix formulation of $V_{n}$ can be given as%
\begin{equation}
\left( 
\begin{array}{c}
V_{n+2} \\ 
V_{n+1} \\ 
V_{n}%
\end{array}%
\right) =\left( 
\begin{array}{ccc}
r & s & t \\ 
1 & 0 & 0 \\ 
0 & 1 & 0%
\end{array}%
\right) ^{n}\left( 
\begin{array}{c}
V_{2} \\ 
V_{1} \\ 
V_{0}%
\end{array}%
\right)  \label{equat:nmouyvbcfxdsz}
\end{equation}%
The matrix $M$ was defined and used in [\ref{bib:shannon1972}]. For matrix
formulation (\ref{equat:nmouyvbcfxdsz}), see [\ref{bib:yalavigi1971}] and [%
\ref{bib:waddill1991}].

Now we define the matrice $M_{V}$ as 
\begin{equation*}
M_{V}=\left( 
\begin{array}{ccc}
\widehat{V}_{4} & \widehat{V}_{3}+\widehat{V}_{2} & \widehat{V}_{3} \\ 
\widehat{V}_{3} & \widehat{V}_{2}+\widehat{V}_{1} & \widehat{V}_{2} \\ 
\widehat{V}_{2} & \widehat{V}_{1}+\widehat{V}_{0} & \widehat{V}_{1}%
\end{array}%
\right) .
\end{equation*}%
This matrice\ $M_{V}$ can be called generalized Tribonacci sedenion matrix.
As special cases, Tribonacci and Tribonacci-Lucas sedenion matricies are%
\begin{equation*}
M_{T}=\left( 
\begin{array}{ccc}
\widehat{T}_{4} & \widehat{T}_{3}+\widehat{T}_{2} & \widehat{T}_{3} \\ 
\widehat{T}_{3} & \widehat{T}_{2}+\widehat{T}_{1} & \widehat{T}_{2} \\ 
\widehat{T}_{2} & \widehat{T}_{1}+\widehat{T}_{0} & \widehat{T}_{1}%
\end{array}%
\right) \text{\ and }M_{K}=\left( 
\begin{array}{ccc}
\widehat{K}_{4} & \widehat{K}_{3}+\widehat{K}_{2} & \widehat{K}_{3} \\ 
\widehat{K}_{3} & \widehat{K}_{2}+\widehat{K}_{1} & \widehat{K}_{2} \\ 
\widehat{K}_{2} & \widehat{K}_{1}+\widehat{K}_{0} & \widehat{K}_{1}%
\end{array}%
\right)
\end{equation*}%
respectively.

\begin{theorem}
For $n\geq 0,$ the following holds:%
\begin{equation}
M_{V}\left( 
\begin{array}{ccc}
r & s & t \\ 
1 & 0 & 0 \\ 
0 & 1 & 0%
\end{array}%
\right) ^{n}=\left( 
\begin{array}{ccc}
\widehat{V}_{n+4} & s\widehat{V}_{n+3}+t\widehat{V}_{n+2} & t\widehat{V}%
_{n+3} \\ 
\widehat{V}_{n+3} & s\widehat{V}_{n+2}+t\widehat{V}_{n+1} & t\widehat{V}%
_{n+2} \\ 
\widehat{V}_{n+2} & s\widehat{V}_{n+1}+t\widehat{V}_{n} & t\widehat{V}_{n+1}%
\end{array}%
\right)  \label{equation:nbhbxczdradweq}
\end{equation}
\end{theorem}

\textit{Proof}. We prove by mathematical induction on $n.$ If $n=0$ then the
result is clear. Now, we assume it is true for $n=k,$ that is%
\begin{equation*}
M_{V}M^{k}=\left( 
\begin{array}{ccc}
\widehat{V}_{k+4} & s\widehat{V}_{k+3}+t\widehat{V}_{k+2} & t\widehat{V}%
_{k+3} \\ 
\widehat{V}_{k+3} & s\widehat{V}_{k+2}+t\widehat{V}_{k+1} & t\widehat{V}%
_{k+2} \\ 
\widehat{V}_{k+2} & s\widehat{V}_{k+1}+t\widehat{V}_{k} & t\widehat{V}_{k+1}%
\end{array}%
\right) .
\end{equation*}%
If we use (\ref{equation:mnouaedcxzqwoup}), then for $k\geq 0,$ we have $%
\widehat{V}_{k+3}=r\widehat{V}_{k+2}+s\widehat{V}_{k+1}+t\widehat{V}_{k}.$
Then by induction hypothesis, we obtain%
\begin{eqnarray*}
M_{V}M^{k+1} &=&(M_{V}M^{k})M \\
&=&\left( 
\begin{array}{ccc}
\widehat{V}_{k+4} & s\widehat{V}_{k+3}+t\widehat{V}_{k+2} & t\widehat{V}%
_{k+3} \\ 
\widehat{V}_{k+3} & s\widehat{V}_{k+2}+t\widehat{V}_{k+1} & t\widehat{V}%
_{k+2} \\ 
\widehat{V}_{k+2} & s\widehat{V}_{k+1}+t\widehat{V}_{k} & t\widehat{V}_{k+1}%
\end{array}%
\right) \left( 
\begin{array}{ccc}
r & s & t \\ 
1 & 0 & 0 \\ 
0 & 1 & 0%
\end{array}%
\right) \\
&=&\left( 
\begin{array}{ccc}
r\widehat{V}_{k+4}+s\widehat{V}_{k+3}+t\widehat{V}_{k+2} & s\widehat{V}%
_{k+4}+t\widehat{V}_{k+3} & t\widehat{V}_{k+4} \\ 
r\widehat{V}_{k+3}+s\widehat{V}_{k+2}+t\widehat{V}_{k+1} & s\widehat{V}%
_{k+3}+t\widehat{V}_{k+2} & t\widehat{V}_{k+3} \\ 
r\widehat{V}_{k+2}+s\widehat{V}_{k+1}+t\widehat{V}_{k} & s\widehat{V}_{k+2}+t%
\widehat{V}_{k+1} & t\widehat{V}_{k+2}%
\end{array}%
\right) \\
&=&\left( 
\begin{array}{ccc}
\widehat{V}_{k+5} & s\widehat{V}_{k+4}+t\widehat{V}_{k+3} & t\widehat{V}%
_{k+4} \\ 
\widehat{V}_{k+4} & s\widehat{V}_{k+3}+t\widehat{V}_{k+2} & t\widehat{V}%
_{k+3} \\ 
\widehat{V}_{k+3} & s\widehat{V}_{k+2}+t\widehat{V}_{k+1} & t\widehat{V}%
_{k+2}%
\end{array}%
\right) .
\end{eqnarray*}%
Thus, (\ref{equation:nbhbxczdradweq}) holds for all non-negative integers $%
n. $%
\endproof%

\begin{corollary}
For $n\geq 0,$ the following holds:%
\begin{equation*}
\widehat{V}_{n+2}=\widehat{V}_{2}U_{n+2}+(\widehat{V}_{1}+\widehat{V}%
_{0})U_{n+1}+\widehat{V}_{1}U_{n}.
\end{equation*}
\end{corollary}

\textit{Proof}. The proof can be seen by the coefficient\ (\ref%
{equation:gbhasmnuopdcxsz}) of the matrix $M_{V}$ and (\ref%
{equation:nbhbxczdradweq}).

Note that we have similar results if we replace the matrix $M$ with the
matrices $N$ and $O$ defined by%
\begin{equation*}
N=\left( 
\begin{array}{ccc}
r & 1 & 0 \\ 
s & 0 & 1 \\ 
t & 0 & 0%
\end{array}%
\right) \text{ and }O=\left( 
\begin{array}{ccc}
0 & 1 & 0 \\ 
0 & 0 & 1 \\ 
r & s & t%
\end{array}%
\right) .
\end{equation*}

\end{document}